\documentstyle[11pt,epic,eepic]{article}
\input epsf.tex
\font\goth=eufm10

\catcode`\@=11
\renewcommand\subsection{\@startsection{subsection}{2}{\z@}%
                                     {-3.25ex\@plus -1ex \@minus -.2ex}%
                                     {-0.01 mm}
                                     {\normalfont\large\bfseries}}
\newcommand{\doublebro}{[\hspace*{-.5mm}[}
\newcommand{\doublebrf}{]\hspace*{-.5mm}]}
\newtheorem{theorem}{Theorem}[section]

\newtheorem{proposition}{Proposition}[section]
\newtheorem{lemma}{Lemma}[section]

\def\boxit#1#2{\setbox1=\hbox{\kern#1{#2}\kern#1}%
\dimen1=\ht1 \advance\dimen1 by #1 \dimen2=\dp1 \advance\dimen2 by #1
\setbox1=\hbox{\vrule height\dimen1 depth\dimen2\box1\vrule}%
\setbox1=\vbox{\hrule\box1\hrule}%
\advance\dimen1 by .4pt \ht1=\dimen1
\advance\dimen2 by .4pt \dp1=\dimen2 \box1\relax}

\def\resp{{\em resp.$\ $}}

\def\cqfd{\hfill $\Box$ \bigskip}
\def\adots{\mathinner{\mkern2mu\raise1pt\hbox{.}
\mkern3mu\raise4pt\hbox{.}\mkern1mu\raise7pt\hbox{.}}}
\def\<{\langle\,}
\def\>{\,\rangle}

\def\cf{{\it cf.$\ $}}

\def\Hom{{\rm Hom\,}}
\def\Ker{{\rm Ker\,}}

\def\ie{{\it i.e. }}

\def\SG{\hbox{\goth S}}

\def\Pad{\hbox{\goth p}}
\def\H{\widehat{H}}

\def\d{{\bf d}}

\def\m{{\bf m}}
\def\kb{{\bf k}}
\def\Nb{{\bf N}}

\def\Hc{{\cal H}}

\def\Fb{{\bf F}}
\def\Fq{{\bf F}_q}

\def\Zb{{\bf Z}}
\def\Cb{{\bf C}}
\def\Cc{{\cal C}}

\def\T2{{\cal T}}

\def\S2{{\cal S}}

\def\mod{{\rm \ mod\ }}
\def\min{{\rm  min }}

\def\g{\hbox{\goth g\hskip 1pt}}

\def\Sl{\hbox{\goth sl\hskip 1pt}}
\def\slchap{\widehat{\hbox{\goth sl}}}

\def\k{{\bf k}}

\def\Stab{{\rm Stab\, }}
\def\Aut{{\rm Aut\, }}
\def\End{{\rm End\, }}
\def\GL{{ GL }}
\def\Mat{{\rm Mat }}
\def\res{{\rm res\,}}

\def\In{{\rm in}}
\def\Out{{\rm out}}
\def\Proof{\medskip\noindent {\it Proof --- \ }}
\def\blambda{{\hbox{\boldmath$\lambda$}}}
\def\deg{{\rm deg\,}}
\def\today{\number\day \space\ifcase\month\or
 Janvier\or F\'evrier\or Mars\or Avril\or Mai\or Juin\or
 Juillet\or Ao\^ut\or Septembre\or Octobre\or Novembre\or D\'ecembre\fi
 \space\number\year}

\title{\bf  Zelevinsky's involution at roots of unity}
\author{\rm Bernard {\sc Leclerc}\thanks{
D\'epartement de Math\'ematiques,
Universit\'e de Caen,
BP 5186, 14032 Caen cedex, France.},
Jean-Yves {\sc Thibon}\thanks{
Institut Gaspard Monge,
Universit\'e de Marne-la-Vall\'ee,
Champs-sur-Marne,
77454 Marne-la-Vall\'ee cedex 2, France.},
\rm and Eric {\sc Vasserot}\thanks{
D\'epartement de Math\'ematiques,
Universit\'e de Cergy-Pontoise, 
2 Av. A. Chauvin, 95302 Cergy-Pontoise
cedex, France.}
}

\begin{document}
\maketitle

%%%%%%%%%%%%%%%%%%%%%%%%%%%%%%%%%%%%%%%%%%%%%%%%%%%%%%%%%%%%%%%%%%%
%  ABSTRACT
%%%%%%%%%%%%%%%%%%%%%%%%%%%%%%%%%%%%%%%%%%%%%%%%%%%%%%%%%%%%%%%%%%%
\vskip 1cm

\begin{abstract}
We give a combinatorial algorithm for computing Zelevinsky's
involution of the set of isomorphism classes of irreducible 
representations of the affine
Hecke algebra $\H_m(t)$ when $t$ is a primitive $n$th root of 1.
We show that the same map can also be interpreted in terms
of aperiodic nilpotent orbits of $\Zb/n\Zb$-graded vector
spaces.
\end{abstract}

\vskip 0.6cm
%%%%%%%%%%%%%%%%%%%%%%%%%%%%%%%%%%%%%%%%%%%%%%%%%%%%%%%%%%%%%%%%%%%
%  SECTION 1
%%%%%%%%%%%%%%%%%%%%%%%%%%%%%%%%%%%%%%%%%%%%%%%%%%%%%%%%%%%%%%%%%%%

\section{Introduction}
In \cite{Zel80}, Zelevinsky has introduced an involution 
of the Grothendieck group of the category of complex smooth
representations of finite length of $G=\GL(m,F)$, where
$F$ is a $\Pad$-adic field, and conjectured that this involution
permutes the classes of irreducible representations
(\cite{Zel80}, 9.17).
In \cite{Zel81}, Zelevinsky further conjectured a geometric description
of this involution in terms of the graded nilpotent orbits that
parametrize the simple $G$-modules.
Both conjectures have been proved by Moeglin and Waldspurger \cite{MW} 
for the category of admissible representations of $G$ generated by
their space of $I$-fixed vectors, where $I$ is an Iwahori subgroup of $G$.
In this case, by a theorem of Bernstein, Borel and Matsumoto, the
conjectures can be reformulated in terms of the Hecke algebra
$\H_m(q)$ of $G$ with respect to $I$. 
This is an affine Hecke algebra of type $A$ with parameter $q$ equal
to the cardinality of the residue field of $F$.
As shown in \cite{MW}, the problem then becomes that of
describing how the isomorphism classes of simple $\H_m(q)$-modules
are permuted when the action is twisted by a certain involutive 
automorphism $\tau$ of $\H_m(q)$. 
In fact the answer does not depend on $q$, and is
the same for any complex parameter $t$ of infinite multiplicative order.

This is no longer the case if one considers the Hecke algebra
$\H_m(t)$ with $t$ a primitive $n$th root of 1. 
According to a conjecture of Vigneras \cite{Vign2}, the solution
of the problem in this case would be relevant for the 
$l$-modular representation theory of $G$, where $l$ is a prime different
from $p$, and $n$ is the multiplicative order of $q$ in the field
with $l$ elements.
 
In this paper, we obtain results similar to those of Moeglin and Waldspurger
in the root of unity case. Namely, (i) given a simple $\H_m(t)$-module $L$,
we describe in Section~5.2 a simple combinatorial algorithm for computing 
the isomorphism class $[L^\tau]$ of the simple module obtained by twisting 
with $\tau$, and (ii) we prove in Section~5.3 that the geometric description  
conjectured by Zelevinsky still holds in this case, this time in terms
of the $\Zb/n\Zb$-graded nilpotent orbits parametrizing the simple
$\H_m(t)$-modules.

However, our methods are different from those of \cite{MW} (and thus we obtain
new proofs of these results in the non root of unity case).
Our main tool is a recent theorem of Ariki \cite{Ar} which relates
the Grothendieck groups of the algebras $\H_m(t)$ with the enveloping
algebra $U^-(\g)$, where $\g$ is the Kac-Moody algebra of type $A_{n-1}^{(1)}$
if $t$ has order $n$, and $A_\infty$ if $t$ has infinite order.
Using this result, we show that our problem is equivalent to describing
the natural 2-fold symmetry of Kashiwara's crystal graph of $U_v^-(\g)$, 
induced by the root diagram automorphism exchanging the simple roots
$\alpha_i$ and $\alpha_{-i}$. 
Then the combinatorial algorithm for $[L^\tau]$ follows immediately
from the explicit description of this graph (Theorem~\ref{TH1}),
and the equivalent formulation in terms of nilpotent orbits is deduced from a
geometric construction of the crystal basis conjectured by Lusztig
\cite{Lu90} and recently verified by Kashiwara and Saito \cite{KS}. 

We note that Procter \cite{Pro} and Aubert \cite{Aub} have established 
independently that the first conjecture of Zelevinsky holds
for the whole category of complex smooth representations of finite 
length of $G$.
Assuming this hypothesis, Moeglin and Walspurger had already shown in 
\cite{MW} that their geometric and combinatorial descriptions 
of Zelevinsky's duality were also 
valid for this larger class of representations. 

The paper is structured as follows. 
In Section~\ref{SECT2} we recall the
definition of Zelevinsky's involution, and the 
parametrization of simple $\H_m(t)$-modules by multisegments.
Then we can formulate our problem in a more precise way.
In Section~\ref{SECT3} we review the Hall-Ringel algebra $\Hc$ associated
with a quiver and the canonical basis of its subalgebra $\Cc$.
In Section~\ref{SECT4}, using the isomorphism $\Cc \cong U_v^-(\g)$,
we compute the action of the Chevalley generators of $U_v^-(\g)$
on the PBW-basis of $\Hc$.
This allows us to describe the crystal graph of $U_v^-(\g)$ in
Lusztig's parametrization by nilpotent orbits, and to derive 
the combinatorial description of $\tau$ (5.2).
Finally, in 5.3 we briefly recall the results of Lusztig, Kashiwara and Saito 
on the construction of the crystal basis in terms of irreducible
components of a certain Lagrangian variety, and we
deduce from them the geometric description of $\tau$.

%%%%%%%%%%%%%%%%%%%%%%%%%%%%%%%%%%%%%%%%%%%%%%%%%%%%%%%%%%%%%
%    SECTION 2
%%%%%%%%%%%%%%%%%%%%%%%%%%%%%%%%%%%%%%%%%%%%%%%%%%%%%%%%%%%%%

\section{Affine Hecke algebras} \label{SECT2}

\subsection{} Let $t$ be a nonzero complex number.
The affine Hecke algebra $\H_m(t)$ associated to $\GL(m)$ is the
associative $\Cb$-algebra with invertible generators
$T_1,\ldots,T_{m-1}$ and $y_1,\ldots,y_m$ subject to the relations
\begin{eqnarray*}
&&T_iT_{i+1}T_i=T_{i+1}T_iT_{i+1},\quad 1\le i\le m-2,\label{EQ_T1}\\
&&T_iT_j=T_jT_i,\qquad \vert i-j\vert>1,\label{EQ_T2}\\
&&(T_i-t)(T_i+1)=0,\quad 1\le i\le m-1,\label{EQ_T3}\\
&&y_iy_j=y_jy_i, \qquad 1\le i,j\le m,\label{EQ_YY}\\
&&y_jT_i=T_iy_j,\qquad \mbox{for $j\not= i,i+1$},\label{EQ_TY2}\\
&&T_iy_iT_i=t\,y_{i+1}, \qquad 1\le i\le m-1 .\label{EQ_Y}
\end{eqnarray*}
The algebra $\H_m(t)$ admits involutive automorphisms 
$\tau$, $\flat$, $\sharp$, given on the generators by
\begin{eqnarray*}
& T_i^\tau = -tT_{m-i}^{-1}, \qquad &y_j^\tau=y_{m+1-j},\label{EQ_TAU}\\
& T_i^\flat = T_{m-i},\qquad & y_j^\flat=y_{m+1-j}^{-1},\label{EQ_bemol}\\
& T_i^\sharp = -tT_i^{-1},\qquad & y_j^\sharp=y_j^{-1}\,.\label{EQ_diese}
\end{eqnarray*}
One checks that $\tau$ and $\flat$ commute and that
${x}^\sharp=(x^\flat)^\tau$.
The involution $\sharp$ was defined by Iwahori and Matsumoto
in \cite{IM}, p. 280, and
$\tau$ is usually called Zelevinsky's involution (see \cite{MW}, p. 146),
because of its relation with the involution 
on representations of $\Pad$-adic $\GL(m)$ introduced by Zelevinsky in 
\cite{Zel80,Zel81}.

\subsection{}
For $\mu=(\mu_1,\ldots,\mu_r)$  a composition of $m$,  set
$D(\mu)=\{\mu_1+\cdots+\mu_k,\ 1\le k\le r-1\}=\{d_1,\ldots,d_{r-1}\}$
and denote
by $\H_\mu$ the subalgebra of $\H_m(t)$ generated by $y_1,\ldots,y_m$ and
$\{T_i\ |\ i\not\in D(\mu)\}$.
For $a\in\Zb^r$, let $\Cb_{\mu,a}$ be the 1-dimensional representation
of $\H_\mu$ defined by $T_i\mapsto t$, $i\not\in D(\mu)$, and 
$y_{d_{i-1}+1}\mapsto t^{a_i}$,
$i=1,\ldots,r$, where $d_0=0$.
We denote by $M_{\mu,a}$ the induced $\H_m(t)$-module
\[
M_{\mu,a} = \H_m(t) \otimes_{\H_\mu} \Cb_{\mu,a} \,.
\]
Let $R(\H_m(t))$ be the complexified Grothendieck group of the category of
finite dimensional $\H_m(t)$-modules, and let $R_m(t)$ be the linear
span of the classes of the composition factors of the various
$M_{\mu,a}$ in $R(\H_m(t))$. We set
\[
R(t)=\bigoplus_{m\ge 0}R_m(t) \,,
\]
where we have put for convenience $R_0(t) = \Cb$.
When $(\nu,b)$ is a permutation of $(\mu,a)$, \ie $\nu_i = \mu_{\sigma(i)}$
and $b_i=a_{\sigma(i)}$ for some $\sigma \in \SG_r$, the induced modules
$M_{\mu,a}$ and $M_{\nu,b}$ are in general non isomorphic, but
their classes in $R(t)$ are equal \cite{Zel80} (see \cite{Rog}, p. 455). 
In the sequel, we shall always work in $R(t)$ and therefore shall
only consider pairs $(\mu,a)$ up to permutation.
In fact, such unordered pairs $(\mu,a)$ are naturally identified with
certain graded nilpotent orbits $O$ (see below Section~\ref{SECT3}),
which gives a canonical labelling of the classes of the induced modules
by these orbits : $[M_{\mu,a}] = [M_O]$.
There is a partial order on orbits given by
$O \unlhd P$ if $O \subset \overline{P}$.

When $t$ is not a root of unity, it is known \cite{Zel80}
that for all $O$, there is a unique simple module whose
class occurs in the expansion of $[M_O]$ but does not occur
in any $[M_P]$ with $O\lhd P$.  
Let $L_O$ denote this simple module. Then the $L_O$ are pairwise
non isomorphic, and one has in $R(t)$
\[
[M_O] = [L_O] + \sum_{O\lhd P} K_{O,P} \, [L_P]
\]
for some positive integers $K_{O,P}$. (The $K_{O,P}$ are in fact values
at 1 of Kazhdan-Lusztig polynomials of type $A$, as conjectured by Zelevinsky 
\cite{Zel81,Zel85} and proved by Ginzburg \cite{CG}, Theorem 8.6.23.)

If $t$ is a primitive $n$th root of unity, we may consider that in a pair
$(\mu,a)$, $a$ belongs to $(\Zb/n\Zb)^r$. Let us say that  $(\mu,a)$
is aperiodic if for each $\ell\in \Nb^*$, the set 
$\{a_i\in \Zb/n\Zb\,|\, \mu_i=\ell\}$
has at most $n-1$ elements. These are the labels associated to the so-called
aperiodic graded nilpotent orbits (see \cite{Lu91}, 15.3). 
It follows from Ariki's theorem \cite{Ar} that for all aperiodic $O$,
there is a unique simple module $L_O$ whose
class occurs in the expansion of $[M_O]$ but does not occur
in any $[M_P]$ with $O\lhd P$.
Moreover, $\{[L_O] \ | \ O \mbox{ aperiodic }\}$ is a basis of $R(t)$.
Finally, the multiplicities of the $[L_P]$ in the $[M_O]$ are
given by certain Kazhdan-Lusztig polynomials of affine type $\tilde A$.
(For the general type a similar result was previously announced by
Grojnowski \cite{Groj}).

\subsection{}
When $t$ is not a root of unity, the pairs $(\mu,a)$ can be 
identified with Zelevinsky's multisegments \cite{Zel80}.
We recall that a segment is an interval $[i,j]$ in $\Zb$, and that
a multisegment is a formal finite unordered sum 
$\m = \sum_{i\le j} m_{ij}[i,j]$.
(Here $m_{ij}$ stands for the multiplicity of the segment $[i,j]$
in $\m$.)
The multisegment corresponding to $(\mu,a)$ is then
$$
\m=\sum_{i=1}^r [a_i,a_i+\mu_i-1]\,.
$$
A multisegment can be conveniently regarded as a coloured
multipartition, \ie as a sequence $(\lambda^{(i)})_{i\in\Zb}$ of
partitions, such that the parts of $\lambda^{(i)}$ are
the lengths of the segments $[i,j]$, and each cell of the $k$th
column of the Young
diagram of $\lambda^{(i)}$ contains the integer $i+k-1$.

\begin{figure}[t]
\begin{center}
\leavevmode
\epsfxsize =10cm
\epsffile{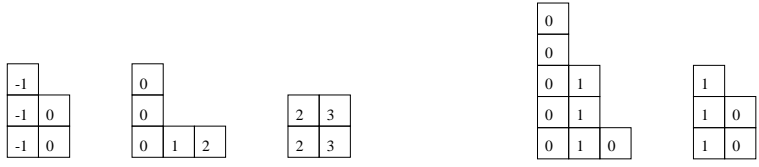}
\end{center}
\caption{\label{FIG1} \small The diagrams of 
$(\mu,a) = ((2,2,3,1,1,2,2,1),(2,2,0,0,0,-1,-1,-1))$ and its reduction 
modulo $2$}
\end{figure}
In the case where $t$ is a primitive $n$th root of unity,
we shall also identify the labels $(\mu,a)$ with multisegments,
this time over $\Zb/n\Zb$, where we regard $\Zb/n\Zb$ as a set
of $n$ points on the unit circle.
For $\ell \in \Nb^*$ and $i \in \Zb/n\Zb$
we define the segment of length $\ell$ and
origin $i$ on this circle (or more appropriately the loop) as:
$$
[i;\ell) := [i, i +1, \ldots , i + \ell -1]\,.
$$ 
Sometimes we shall also need the dual notation
$$
(\ell;i] := [i-\ell+1, i-\ell+2, \ldots , i]\,.
$$
A multisegment over $\Zb/n\Zb$ is now a formal finite sum
$$
\m = \sum_{i\in\Zb/n\Zb, \ \ell\in \Nb^*} m_{[i;\ell)} \, [i;\ell)
   = \sum_{j\in\Zb/n\Zb, \ k\in \Nb^*} m_{(k;j]} \, (k;j]\,.
$$
The cyclic multisegment corresponding to the label $(\mu,a)$ is 
$\m = \sum_i [a_i;\mu_i)$.
These cyclic multisegments may also be represented by $n$-tuples
of partitions $(\lambda^{(i)})_{i\in\Zb/n\Zb}$, with cells
labelled by integers modulo $n$.
A label $(\mu,a)$ is aperiodic if for each $\ell$, there is
at least one component $\lambda^{(i)}$ which does not contain
the part $\ell$ (this is, up to conjugation of partitions, the labelling used by
Ringel in \cite{Ri93} for a basis of the composition
algebra of the cyclic quiver). 
These definitions are illustrated in Figure~\ref{FIG1}.

Let $I = \Zb$ (\resp $\Zb/n\Zb$). 
The set of all multisegments over $I$ will be denoted by ${\cal M}(I)$,
or simply ${\cal M}$ if no confusion can arise.
Thus we have a parametrization of the simple $\H_m(t)$-modules $L_\m$
by $\m \in {\cal M}(I)$ (with $\m$ aperiodic for $I=\Zb/n\Zb$).
We introduce a $\Nb^{(I)}$-grading of $R(t)$ by putting
$\deg [L_\m] = (d_i)_{i\in I}$, where $d_i$ is the number of cells
of the diagram of $\m$ which contain $i$.

\subsection{}
The three involutions $\tau$, $\flat$ and $\sharp$ induce
involutions on the set of simple $\H_m(t)$-modules 
$L_\m$, and hence on the set of multisegments $\m$ that parametrize them.
It is easy to see that $L_\m^\flat=L_{\m^\flat}$ where if 
$\m = \sum_{i\in I, \ \ell\in \Nb^*} m_{[i;\ell)} \, [i;\ell)$,
then $\m = \sum_{i\in I, \ \ell\in \Nb^*} m_{[i;\ell)} \, (\ell;-i]$.
Therefore, it is equivalent to describe either $\tau$ or $\sharp$.
For generic $t$, a geometric description of $\tau$ in terms of graded
nilpotent orbits has been conjectured by Zelevinsky \cite{Zel81}.
This conjecture was proved by Moeglin and Waldspurger \cite{MW},
who also gave a combinatorial algorithm for computing $\m \mapsto \m^\tau$.
An explicit formula  for $\m^\tau$ has also been found by Knight and 
Zelevinsky \cite{KZ}.

The finite Hecke algebra $H_m(t)$ is the quotient  of $\H_m(t)$
by the relation $y_1=1$. Hence,  $\sharp$ induces
an involution on simple $H_m(t)$-modules. 
When $t$ is not a root
of unity, these modules are parametrized by partitions of $m$,
and $\sharp$ is just the conjugation of partitions.
When $t$ is an $n$th root of unity, the simple $H_m(t)$-modules are
labelled by $n$-regular partitions, \ie partitions
with no part repeated more than $n-1$ times, and it was recently proved
that $\sharp$ coincides with a bijection previously defined by Mullineux \cite{Mu}.
One proof follows from \cite{LLT96} p. 229 and Ariki's theorem, and 
a second proof was obtained by Brundan \cite{Br} by ``quantizing" the results of
Kleshchev on modular representations of symmetric groups. 

The aim of this paper is to
describe the maps $\tau$ and $\sharp$ in the affine case
when $t$ is a primitive $n$th root of unity.
Proposition~\ref{PRO5.1} below can be regarded as a common
generalization of Mullineux's bijection as described by Kleshchev
\cite{Kl3}, and of the algorithm of Moeglin and Waldspurger.

%%%%%%%%%%%%%%%%%%%%%%%%%%%%%%%%%%%%%%%%%%%%%%%%%%%%%%%%%%%%%%%%
% SECTION 3
%%%%%%%%%%%%%%%%%%%%%%%%%%%%%%%%%%%%%%%%%%%%%%%%%%%%%%%%%%%%%%%%%
 
\section{Graded nilpotent orbits and Hall-Ringel algebras} \label{SECT3}

\subsection{}
Let $\Gamma$ denote the quiver of type $A_{n-1}$ 
(\resp $A_\infty$ or $A_{n-1}^{(1)}$). 
The set of vertices of $\Gamma$ is $I=[1,n-1]$ (\resp 
$\Zb$ or $\Zb/n\Zb$), and the set of arrows is 
$\Omega=\{i\to i-1\ |\ i,i-1\in I\}$. 
A nilpotent representation of $\Gamma$ over the field $\k$ is a pair $(V,x)$,
where $V = \oplus_{i\in I} V_i$ is an $I$-graded finite-dimensional vector space,
and $x=(x_a)_{a\in\Omega}$ is a degree -1 nilpotent linear operator
on $V$. In other words, $x$ is an element of $N_{V,\Omega}$,
the subset of nilpotent elements of
$$
E_{V,\Omega}=\bigoplus_{i\rightarrow j\in\Omega}\Hom(V_i,V_j)\,.
$$
The vector $\d = {\bf dim}\, V =(\dim V_i)_{i\in I}\in\Nb^{(I)}$ 
is the (graded) dimension of the representation.
The group $G_V=\prod_{i\in I}\GL(V_i)$ acts on $E_{V,\Omega}$ and $N_{V,\Omega}$
by conjugation.
Two representations $(V,x)$ and $(V,x')$ 
are equivalent if $x$ and $x'$ lie in the same orbit. 
Thus, the set of isomorphism classes of nilpotent $\Gamma$-representations
is naturally identified with ${\cal O} = \sqcup_{\d \in \Nb^{(I)}} {\cal O_\d}$,
where ${\cal O_\d}$ is the set of $G_V$-orbits in $N_{V,\Omega}$ for
any graded $\k$-vector space $V$ of dimension $\d$. 

The isomorphism classes of nilpotent $\Gamma$-representations
are parametrized by ${\cal M}(I)$ in the following way.
For each vertex $i\in I$, let $\k[i]$ denote the simple
$\Gamma$-module for which $V=V_i=\k$ and $x=0$.
Given a positive integer $\ell$, there is a unique 
(up to isomorphism) indecomposable $\Gamma$-module 
$\k[\ell;i]$ of length $\ell$ with head $\k[i]$,
and all indecomposable modules are of this type.
Therefore, any nilpotent $\Gamma$-module $M$ is isomorphic to 
\[
\k[\m] := \bigoplus_{i\in I, \ell \in \Nb^*}
\k[\ell;i]^{\oplus m_{(\ell;i]}} 
\]
for a unique multisegment $\m$.
We denote by $O_\m$ the $G_V$-orbit of the module $\k[\m]$. 
For simplicity, we also denote by $O_i$ the orbit of the simple module
$\k[i]$.

\subsection{}
We recall the definition of the twisted Hall algebra associated 
with the quiver $\Gamma$ (see \cite{Ri90,Rin93,Lu98}).
The classification of nilpotent $\Gamma$-modules is independent 
of the ground field.
Take $\kb = \Fq$, the field with $q$ elements, and let $O,P,Q$
be orbits in $N_{V,\Omega}, N_{W,\Omega}, N_{U,\Omega}$
respectively, where $\dim U = \dim V + \dim W$.
Let $M$ be a $\Gamma$-module lying on the orbit $Q$.
Then it is was proved by Ringel \cite{Ri90,Ri93}
that the number of submodules of $M$ with type $P$
and cotype $O$ is a polynomial in $q$ with integer coefficients,
independent of the choice of $M$ on $Q$.
This is the Hall polynomial $F_{O,P}^{Q}(q)$.

Define a bilinear form $m$ on $\Nb^{(I)}$ by
\begin{equation}
m({\bf a},{\bf b})=\sum_{i\to j\in\Omega}a_ib_j+\sum_{i\in I}a_ib_i\,.
\end{equation}
Then the generic twisted Hall algebra associated with $\Gamma$ is the 
$\Cb(v)$-algebra $\Hc$ with basis $\{u_O\ | \ O \in {\cal O}\}$
and multiplication
\begin{equation}\label{MULT1}
u_O\circ u_P=v^{m(\dim V,\dim W)}\sum_Q F^Q_{O,P}(v^{-2})\,u_Q.
\end{equation}
(In the setting of \cite{Lu98}, $u_O$ stands for the product
of characteristic functions of the sets of $\Fb_{p^{2m}}$-rational points
of $O$ for a fixed prime $p$.)

\subsection{}
Let $O$ be a $G_V$-orbit in $N_{V,\Omega}$. 
One sets
\begin{equation}
\< O \> = v^{\dim O} \, u_O \,,
\end{equation}
where $\dim O$ denotes the dimension of the orbit $O$
(not to be confused with the dimension $\dim V$ of 
the representation of $\Gamma$ corresponding to a point of $O$).
The basis $\{ \< O \> \}$ is called the Poincar\'e-Birkoff-Witt
basis of $\Hc$.
If $P$ is a $G_W$-orbit in $N_{W,\Omega}$, it follows from (\ref{MULT1}) that
\begin{equation} \label{MULT2}
\< O \>\circ \< P \> =\sum_Q v^{\alpha(O,P,Q)}\, F^Q_{O,P}(v^{-2})\, \< Q \>,
\end{equation}
where
\begin{equation}
\alpha(O,P,Q) = \dim O + \dim P - \dim Q + m(\dim V , \dim W) \,.
\end{equation}
In the sequel, we shall use another expression for $\alpha(O,P,Q)$.
For ${\bf a}, {\bf b} \in \Nb^{(I)}$, set 
\begin{equation}
r({\bf a},{\bf b})=-\sum_{i\to j\in\Omega}a_ib_j+\sum_{i\in I}a_ib_i\,,
\end{equation}
and denote by $\varepsilon(O)$ the dimension of the space of endomorphisms
of a module $M$ lying in $O$.
Then
\begin{equation}\label{TWISTR}
\alpha(O,P,Q) = - \varepsilon(O) -\varepsilon(P) + \varepsilon(Q) 
- r(\dim V , \dim W) \,.
\end{equation}
Indeed, $$\dim O = \dim G_V - \dim (\Stab M)\,,$$ and
$$\dim (\Stab M) = \dim (\Aut M) = \dim (\End M) = \varepsilon(O)\,.$$
Thus, if we set $\dim V = (d_{i,V})_{i\in I}$ 
and $\dim W = (d_{i,W})_{i\in I}$, we have
$$
\alpha(O,P,Q) =  - \varepsilon(O) -\varepsilon(P) + \varepsilon(Q)
+ \sum_i d_{i,V}^2 + \sum_i d_{i,W}^2 -  \sum_i (d_{i,V} + d_{i,W})^2
+ m(\dim V , \dim W)
$$
$$
= - \varepsilon(O) -\varepsilon(P) + \varepsilon(Q)
- r(\dim V , \dim W) \,.
$$

\subsection{}
Let ${\cal C}$ be the twisted composition algebra of $\Gamma$, {\it i.e.}
the subalgebra of ${\cal H}$ generated by the characteristic
functions $u_i=\<O_i\>$ of the orbits of the simple modules
$\k[i]$. When $\Gamma$ is of type $A_{n-1}$ or $A_\infty$,
one has $\Cc = \Hc$, but 
$\Cc$ is strictly contained in $\Hc$ for the type $A_{n-1}^{(1)}$. 
The vector space $\Hc$ has a basis $\{b_O \ | \ O \in {\cal O}\}$ 
given by
\begin{equation}\label{CANO}
 b_O = \sum_{i, O'}
v^{-i +\dim O-\dim O'}
\dim {\cal H}^i_{O'}\left( IC_{O} \right)
\<O'\> \,
\end{equation}
where ${\cal H}^i_{O'}\left( IC_O \right)$
is the stalk at a point of $O'$ of the $i$th intersection cohomology
sheaf of the closure $\overline{O}$ of $O$.

For $\Gamma$ of type $A_\infty$, these sheaves have been first 
considered by Zelevinsky in \cite{Zel81}, where he made a conjecture
about their connection with representations of $\Pad$-adic $GL(m)$
known as the $\Pad$-adic analogue of the Kazhdan-Lusztig conjecture.
In \cite{Zel85}, Zelevinsky further proved that the Poincar\'e
polynomials 
$\sum_i t^i \dim {\cal H}^i_{O'}\left( IC_{O} \right)$
are equal to certain Kazhdan-Lusztig polynomials 
$P_{w(O')w(O)}(t^2)$ of type $A$.

For $\Gamma$ of type $A_{n-1}^{(1)}$, Lusztig has shown 
that the Poincar\'e polynomials are also Kazhdan-Lusztig
polynomials, this time of affine type $\tilde A$
(\cite{Lus90},11).
In \cite{Lu91}, he introduced  
the subfamily of aperiodic nilpotent orbits. 
These are the orbits $O_\m$ parametrized by aperiodic multisegments
(see 2.3).  
It is proved in \cite{Lu92}, Section 5 that $\{b_O \ | \ O\ \mbox{aperiodic}\}$ 
is a basis of $\Cc$. 
Note that non aperiodic orbits $O'$ may occur in the right-hand side of 
(\ref{CANO}) for an aperiodic orbit $O$.

Let $x\mapsto \bar x$ be the  ring automorphism of
$\Cc$ defined by $\bar u_i=u_i$ and $\bar v = v^{-1}$.
Let ${\cal L}$ be the $\Zb[v]$-submodule of $\Hc$
spanned by the PBW basis $\<O_{\bf m}\>$. 
Let $O \in {\cal O}$ and suppose that $O$ is aperiodic if $\Gamma$
is of type $A_{n-1}^{(1)}$.
Then $b_O$
is characterized as the unique element of ${\cal L}\cap \Cc$ such
that $\overline{ b_O} = b_O$ and $b_O\equiv \<O\> \mod v{\cal L}$.

%%%%%%%%%%%%%%%%%%%%%%%%%%%%%%%%%%%%%%%%%%%%%%%%%%%%%%%%%%%%%%%%%
% SECTION 4
%%%%%%%%%%%%%%%%%%%%%%%%%%%%%%%%%%%%%%%%%%%%%%%%%%%%%%%%%%%%%%%%%

\section{Affine Lie algebras and quantum affine algebras} \label{SECT4}

\vskip3mm

\subsection{}
Let $\g$ be the Kac-Moody Lie algebra associated with the quiver
$\Gamma$ of type $A_{n-1}$ (\resp $A_\infty$ or $A_{n-1}^{(1)}$),
namely, $\g = \Sl_n$ (\resp $\g = \Sl_\infty$ or $\slchap_n$).
Let $U_v(\g)$ be the corresponding quantized universal enveloping
algebra, with generators $e_i, f_i, v^{h_i}, (i\in I)$.
The subalgebra generated by $f_i, (i\in I)$ is denoted 
by $U_v^-(\g)$.
This algebra is isomorphic to the twisted composition algebra $\Cc$ of $\Gamma$,
the isomorphism being given by $f_i \mapsto \< O_i \>$.
It follows that the twisted Hall algebra $\Hc$ can be regarded as a left
$U_v^-(\g)$-module.
The basis $\{b_O\}$ of $\Cc$ becomes via this isomorphism
a basis of $U_v^-(\g)$ : this is Lusztig's canonical basis.
It is known that it coincides with Kashiwara's lower global basis
of $U_v^-(\g)$.
By taking $v=1$ one obtains the canonical basis of $U^-(\g)$,
that we shall still denote by $\{b_O\}$.
The aim of this section is to describe Kashiwara's crystal graph of $\Hc$.
To do so, we shall first obtain explicit formulas for the action
of the Chevalley generators $f_i$ and their adjoint $e'_i$ 
on $\Hc$.
We believe that these formulas are of independent interest.
For example they yield an algorithm similar to that of \cite{LLT96}
for computing the canonical basis in its expansion on the PBW-basis.

\subsection{}
We describe the action of the $f_i$ on $\Hc$ in the PBW-basis.
We assume from now on that $\Gamma$ is of type $A_{n-1}^{(1)}$.
(The other two cases are readily deduced from this one, since 
$U_v(\Sl_n)$ is a natural subalgebra of $U_v(\slchap_n)$,
and $U_v(\Sl_\infty)$ is the ``limit $n\rightarrow \infty$ of $U_v(\Sl_n)$").
 
As explained above, the PBW-basis is naturally labelled by
multisegments $\m$ over $\Zb/n\Zb$.
For simplicity, we abuse notation and write $\< \m \> := \< O_\m \>$
for the corresponding element of the PBW-basis.
Given $i\in \Zb/n\Zb$, $\ell \in \Nb^*$, and a multisegment $\m$ 
such that $m_{(\ell-1;i-1]}\not = 0$ if $\ell > 1$, we define a new multisegment
$$
\m^+_{\,\ell,i} = \left| \matrix{
&\m + (1;i] &\hbox{ if $\ell=1$,} \cr
&\m + (\ell;i] - (\ell-1;i-1] &\hbox{ if $\ell>1$.}
}\right.
$$
In case $m_{(\ell-1;i-1]} = 0$, we put $\<\m^+_{\,\ell,i}\> = 0$.
\begin{proposition}\label{PRO4.1}
The Chevalley generators of $U_v^-(\slchap_n)$ act on the PBW-basis
of $\Hc$ by
$$
f_i\, \<\m\> = \sum_{\ell\in \Nb^*} 
v^{\,\sum_{k>\ell}(m_{(k-1;i-1]} - m_{(k;i]})}
\,[m_{(\ell;i]} + 1] \, \<\m^+_{\,\ell,i}\> \,,
$$
where for $a \in \Zb$, $[a] = (v^a-v^{-a})/(v-v^{-1})$ is the usual 
$v$-integer.
\end{proposition}
\Proof
We have to apply formulas (\ref{MULT2}) and (\ref{TWISTR}) 
with $O = O_i$ and $P=O_\m$.
We see easily that here, $\varepsilon(O) = 1$, and
$$
r(\dim V , \dim W) = \sum_{k\in \Nb^*} m_{[i;k)} 
- \sum_{k\in \Nb^*} m_{(k;i-1]} \,.
$$
Also, it is clear that in this case, the orbits $Q$ giving a nonzero contribution 
to the right-hand side of (\ref{MULT2}) are those parametrized by
the multisegments of the form $\m^+_{\,\ell,i}$ for some $\ell\in \Nb^*$.
%Let $\kb[\m]$ denote the $\Gamma$-module associated
%to a multisegment $\m$.
We distinguish two cases depending on whether $\ell=1$ or $\ell>1$. 

In the first case we have
\begin{eqnarray*}
\varepsilon(Q) - \varepsilon(P)
& = & \dim (\End (\kb[i])) + \dim(\Hom (\kb[i],\kb[\m] ))
+ \dim (\Hom (\kb[\m],\kb[i]))
\\[2mm]
&=& 1 + \sum_{k\in \Nb^*} m_{[i;k)}
+ \sum_{k\in \Nb^*} m_{(k;i]}\,.
\end{eqnarray*}
To calculate the Hall polynomial $F_{O, P}^Q$,
let us count the number of submodules $(V,x)$ of $(V',x')=\Fq[\m^+_{\,1,i}]$
isomorphic to $\Fq[\m]$.
For $j \not = i$ we must have $V_j = V'_j$, and 
$V_i$ must be a hyperplane in $V'_i$ 
which (i) contains $x(V'_{i+1})$ (so that $V$ be stable under $x'$)
and (ii) does not contain $\Ker x'|_{V'_i}$
(so that $(V,x)$ be isomorphic to $\Fq[\m]$).
The number of hyperplanes satisfying (i) is clearly 
$\doublebro 1 + \sum_{k\in\Nb^*} m_{(k;i]} \doublebrf$
and the number of those satisfying (i) but not (ii) is
$\doublebro \sum_{k\ge 2} m_{(k;i]} \doublebrf$.
(Here we have used the notation 
$\doublebro a \doublebrf = (q^a-1)/(q-1)$.)
Therefore,
$$
F_{O,P}^Q(q) =
\doublebro 1 + \sum_{k\in\Nb^*} m_{(k;i]} \doublebrf
-
\doublebro \sum_{k\ge 2} m_{(k;i]} \doublebrf
=
q^{\sum_{k\ge 2} m_{(k;i]}}\, \doublebro 1 + m_{(1;i]} \doublebrf \,,
$$
and finally, the coefficient of $\<\m^+_{\,1,i}\>$ in $f_i\, \<\m\>$
is equal to 
$$
v^{-1 + 1 +\sum_{k\in \Nb^*} m_{[i;k)} + \sum_{k\in \Nb^*} m_{(k;i]}
-\sum_{k\in \Nb^*} m_{[i;k)} + \sum_{k\in \Nb^*} m_{(k;i-1]}
-2\sum_{k\ge 2} m_{(k;i]}
-m_{(1;i]}}
[m_{(1;i]} + 1 ] 
$$
$$
= \ v^{\,\sum_{k>1}(m_{(k-1;i-1]} - m_{(k;i]})}
\,[m_{(1,i]} + 1]\,.
$$

In the case $\ell > 1$, the calculation is similar and one gets
\begin{eqnarray*}
\varepsilon(Q) - \varepsilon(P)
& = & 1 + \sum_{k\in \Nb^*} m_{[i;k)}
+ \sum_{k\ge \ell} m_{(k;i]}
- \sum_{k\le \ell-1} m_{(k;i-1]}\,, \\[2mm]
F_{O, P}^Q(q) & = &
\doublebro 1 + \sum_{k\ge \ell} m_{(k;i]} \doublebrf
-
\doublebro \sum_{k\ge \ell +1} m_{(k;i]} \doublebrf %\\[2mm]
=
q^{\sum_{k\ge \ell +1} m_{(k;i]}}\, \doublebro 1 + m_{(\ell;i]} \doublebrf \,,
\end{eqnarray*}
so that the coefficient of $\<\m^+_{\,\ell,i}\>$ in $f_i\, \<\m\>$
is equal to
$$
v^{-1 + 1 +\sum_{k\in \Nb^*} m_{[i;k)} + \sum_{k\ge \ell} m_{(k;i]} - 
\sum_{k\le \ell-1} m_{(k;i-1]}
-\sum_{k\in \Nb^*} m_{[i;k)} + \sum_{k\in \Nb^*} m_{(k;i-1]}}
\qquad \qquad
$$
$$
\qquad \qquad \qquad \qquad \qquad \qquad \qquad \qquad
\times \ 
v^{-2\sum_{k\ge \ell+1} m_{(k;i]}}\,
v^{-m_{(\ell;i]}}\,
[m_{(\ell;i]} + 1 ]
$$
$$
= \ v^{\,\sum_{k>\ell}(m_{(k-1;i-1]} - m_{(k;i]})}
\,[m_{(\ell,i]} + 1]\,.
$$
\cqfd

\bigskip
\noindent
{\bf Remark} \quad The specialization $v=1$ of Proposition~\ref{PRO4.1} 
appears in \cite{Ar} Lemma~4.2.

\subsection{}
To define the Kashiwara operators $\tilde f_i$
and $\tilde e_i$ of the $U_v^-(\slchap_n)$-module $\Hc$,
one needs to introduce endomorphisms $e'_i$ of the space $\Hc$
(see \cite{Ka}, 3.4).
First one defines a scalar product on $\Hc$ for which the basis
$\{\<\m\>\}$ is orthogonal by 
\begin{equation}\label{EQ22}
(\<\m \> \,,\, \< {\m'} \>) =
v^{-2\varepsilon(\m)} \, {(1-v^2)^{\dim \k[\m ]}
\over
\left |\Aut \Fq [\m] \right |_{q=v^{-2}} }
\, \delta_{\m\,,\, \m'} \,.
\end{equation}
Here we have written for short $\varepsilon(\m) = \varepsilon(O_\m)$.
This scalar product is the one introduced by Green \cite{Gre},
multiplied by a normalization factor $(v^{-2}-1)^{\dim \k[\m ]}$ 
so that it induces Kashiwara's scalar
product \cite{Ka} on $\Cc \cong U_v^-(\slchap_n)$. 

\begin{lemma} \label{LEM4.1} Put $\varphi_m(t) = (1-t^m)(1-t^{m-1})\cdots (1-t)$.
We have
$$
\left |\Aut \Fq [\m] \right | = \left |\End \Fq [\m] \right | \, 
\prod_{i,\ell} 
{\left |\GL(m_{(\ell ; i]},\Fq) \right |
\over
\left |\Mat(m_{(\ell ; i]},\Fq) \right |}
= q^{\varepsilon(\m)}\, \prod_{i,\ell} \varphi_{m_{(\ell;i]}}(q^{-1})
\,.
$$
\end{lemma}
\Proof
Let $(V,x)$ be a $\Gamma$-module isomorphic to $\Fq[\m]$, and
let $W_{\ell,i}$ denote the subspace of $V_i$ consisting of
all the generators of an indecomposable summand of $(V,x)$ of type $(\ell;i]$.
Thus, $\dim W_{\ell,i} = m_{(\ell;i]}$.

An endomorphism $\varphi$ of $(V,x)$ is completely determined 
by its restriction $\varphi_{|W}$ to $W = \oplus_{\ell,i} W_{\ell,i}$.
Let $U_i$ be a complement of $\oplus_\ell W_{\ell,i}$ in $V_i$, 
and $p_{\ell,i}$ be the corresponding projection of $V_i$
onto $W_{\ell,i}$. Then, $x \circ \varphi = \varphi \circ x$ implies
that $p_{\ell,i} \circ \varphi_{|W_{k,i}} = 0$ for $\ell > k$.
On the other hand, $p_{\ell,i} \circ \varphi_{|W_{\ell,i}}$ is not 
submitted to any condition and may be an arbitrary element
of $\End W_{\ell,i} \cong \Mat(m_{(\ell;i]},\Fq)$.

Now, $\varphi$ is an automorphism of $(V,x)$ if and only if 
$p_W \circ \varphi_{|W}$ is in $\GL(W)$, where $p_W=\sum_{\ell,i} p_{\ell,i}$.
Hence $\varphi$ is an automorphism if and only if
$p_{\ell,i} \circ \varphi_{|W_{\ell,i}}$ belongs to
$\GL(W_{\ell,i}) \cong \GL(m_{(\ell;i]},\Fq)$, and the lemma follows. \cqfd

Using Lemma~\ref{LEM4.1} and the formula
$\dim \k[\m ] = \sum_{i,\ell} \ell \, m_{(\ell ;i]} $, we rewrite (\ref{EQ22})
as
\begin{equation}\label{SP}
(\< \m \> \,,\, \< \m' \>) =
\prod_{i,\ell}
v^{-{m_{(\ell;i]} \choose 2}} \, {(1-v^2)^{(\ell-1)m_{(\ell;i]}}
\over
[m_{(\ell;i]}]! }
\, \delta_{\m\,,\, \m'} \,.
\end{equation}

Next we define $e'_i \in \End_{\Cb(v)} \Hc$ as the adjoint operator
of $f_i$ with respect to this scalar product.
The following proposition describes the action of $e'_i$ on the PBW-basis.
Given $i\in \Zb/n\Zb$, $\ell \in \Nb^*$, and a multisegment $\m$
such that $m_{(\ell;i]}\not = 0$, we put 
$$
\m^-_{\,\ell,i} = \left| \matrix{
&\m - (1;i] &\hbox{ if $\ell=1$,} \cr
&\m - (\ell;i] + (\ell-1;i-1] &\hbox{ if $\ell>1$.}
}\right.
$$
In case $m_{(\ell;i]} = 0$, we put $\<\m^-_{\,\ell,i}\> = 0$.

\begin{proposition}\label{PRO4.2}
The endomorphisms $e'_i$ act on the PBW-basis of $\Hc$ by
$$
e'_i\, \<\m\> = 
v^{\,\sum_{k>1}(m_{(k-1;i-1]} - m_{(k;i]})}
\,v^{-m_{(1;i]} +1}
\,\<\m^-_{\,1,i}\>
\qquad \qquad \qquad \qquad \qquad \qquad 
$$
$$
\qquad \qquad \qquad \qquad  
+ \ 
\sum_{\ell \ge 2}
v^{\,\sum_{k>\ell}(m_{(k-1;i-1]} - m_{(k;i]})}
\,v^{-m_{(\ell;i]} +1}
\,(1-v^{2(m_{(\ell-1;i-1]}+1)})
\,\<\m^-_{\,\ell,i}\> \,.
$$
\end{proposition}
\Proof
By definition of $e'_i$ and using Proposition~\ref{PRO4.1}, 
the coefficient of $\<\m^-_{\,\ell,i}\>$
in $e'_i \, \<\m\>$ is equal to 
\[
{(\< \m \> , f_i \<\m^-_{\,\ell,i}\> )
\over
(\<\m^-_{\,\ell,i}\> , \<\m^-_{\,\ell,i}\> )}
=
{(\< \m \> , \<\m\> )
\over
(\<\m^-_{\,\ell,i}\> , \<\m^-_{\,\ell,i}\> )}
\,
v^{\,\sum_{k>\ell}(m_{(k-1;i-1]} - m_{(k;i]})}
\,[m_{(\ell;i]}]\,.
\]
Now by (\ref{SP}), if $\ell > 1$
\[
{(\< \m \> , \<\m\> )
\over
(\<\m^-_{\,\ell,i}\> , \<\m^-_{\,\ell,i}\> )}
=
(1-v^2)\,v^{-m_{(\ell;i]} + m_{(\ell-1;i-1]} + 1} \,
{[m_{(\ell-1;i-1]} + 1] \over [m_{(\ell;i]}]}
\,,
\]
and if $\ell = 1$
\[
{(\< \m \> , \<\m\> )
\over
(\<\m^-_{\,1,i}\> , \<\m^-_{\,1,i}\> )}
=
{v^{-m_{(1;i]} + 1} \,
\over [m_{(1;i]}]}
\,,
\]
which gives the required result. \cqfd

\bigskip
\noindent
{\bf Remark} \quad  Propositions \ref{PRO4.1}, \ref{PRO4.2}
should be compared to the formulas of \cite{JMMO} for the
action of the Chevalley generators of $U_v(\slchap_n)$
on the level $l$ Fock spaces. Actually, our formulas can be 
regarded as the limit $l \rightarrow \infty$ of those
of \cite{JMMO}.

\subsection{}
Let $\widehat{B}(\infty)$ be Kashiwara's crystal graph of
$\Hc$ for type $A_{n-1}^{(1)}$ in Lusztig's geometric parametrization.
That is, the vertices of this graph are the multisegments $\m$
(or the corresponding orbits $O_\m$),
and there is an arrow $\m {\buildrel i\over\longrightarrow} \m'$
in the graph
if and only if $\tilde f_i  \<O_\m \> \equiv \<O_{\m'}\> \mod v{\cal L}$,
where $\tilde f_i$ denotes as usual the Kashiwara operator \cite{Ka}.
In particular, the connected component $B(\infty)$
of this graph containing the empty multisegment gives the
crystal graph of $\Cc \cong U_v^-(\slchap_n)$ parametrized by the set
of aperiodic orbits.

There are several descriptions of the crystal graph of $U_v^-(\g)$
available in the literature (see \cite{Ka93}, \cite{KKM}, \cite{NZ97}), 
but the only ones using Lusztig's parametrization are given by Reineke 
\cite{Re97} and Kashiwara-Saito \cite{KS}.
Unfortunately, Reineke's description is restricted to finite dimensional $\g$,
while Kashiwara-Saito's description is in terms of the geometry
of orbits, and does not seem to yield a combinatorial algorithm
(see Section 5.3).
Our description of $\widehat{B}(\infty)$ will be an extension to type 
$A_{n-1}^{(1)}$ of Reineke's result for type $A_{n-1}$.

\begin{theorem}\label{TH1}
Let $\m$ be a multisegment over $\Zb/n\Zb$, and fix $i\in \Zb/n\Zb$.
For $k\in \Nb$, put 
$S_{k,i} = \sum_{\ell \ge k} (m_{(\ell;i-1]} - m_{(\ell;i]})$
and let $k_0$ be minimal such that $S_{k_0,i}=\min_k\, S_{k,i}$.
Then there is an arrow $\m {\buildrel i\over\longrightarrow} \m^+_{\,k_0,i}$
in the crystal graph $\widehat{B}(\infty)$ of $\Hc$.
\end{theorem}

\Proof
Our strategy is to reduce the problem to Reineke's description of the crystal
graph of $U^-_v(\Sl_n)$.
We first note that one can let $n$ tend to infinity in Reineke's work
and get the crystal graph of $U^-_v(\Sl_\infty)$.

Let us identify $i\in \Zb/n\Zb$ with its representative in $\{0,1,\ldots ,n-1\}$.
Choose $n$ integers $j_0,j_1,\ldots, j_{n-1}$ such that
$j_r \equiv r \mod n$
and
$j_{i-1} = i-1, \ j_i = i$.
For example  $\{j_r\} = \{0,1,\ldots , n-1\}$ if $i\not =0$,
and $\{j_r\} = \{-1,0,1,\ldots , n-2\}$ if $i=0$.
Define an embedding $\phi_i$ of the set of $\Zb/n\Zb$-multisegments
in the set of $\Zb$-multisegments by
\[
\phi_i\left(\sum_{\ell \in \Nb^*}\sum_{r\in \Zb/n\Zb} m_{(\ell;r]}\,(\ell;r]\right)
=
\sum_{\ell \in \Nb^*}\sum_{r\in \Zb/n\Zb}  m_{(\ell;r]}\,(\ell;j_r]\,.
\]
Then define a $\Cb(v)$-linear embedding
$\Phi_i \ : \  \Hc \longrightarrow U^-_v(\Sl_\infty)$ by
\[
\Phi_i(\<O_\m\>) = \<O_{\phi_i(\m)}\> \,.
\]
Let $f_i^\infty$ denote the Chevalley element in $U^-_v(\Sl_\infty)$
and $e_i^{' \infty}$ its adjoint. Proposition~\ref{PRO4.1} and \ref{PRO4.2}
readily imply  
\begin{lemma} \label{LEMG}
{\rm (i)} The subspace $\Phi_i(\Hc)$ of $U^-_v(\Sl_\infty)$
is stable under $f_i^\infty$ and $e_i^{' \infty}$.
\hfill\break
{\rm (ii)} $\quad \Phi_k \circ f_i = f_i^\infty \circ \Phi_i\,, \qquad
\Phi_i \circ e_i' = e_i^{' \infty} \circ \Phi_i\,.$
\end{lemma}
Let $\gamma_i$ (\resp $\gamma_i^{\infty}$)
be the subgraph of the crystal of $\Hc$ (\resp $U^-_v(\Sl_\infty)$)
obtained by erasing all arrows labelled by $j\not = i$.
Lemma~\ref{LEMG} shows that 
\[
\quad \Phi_i \circ \tilde f_i = \tilde f_i^\infty \circ \Phi_i\,, \qquad
\Phi_i \circ \tilde e_i = \tilde e_i^\infty \circ \Phi_i\,.
\]
Thus, $\phi_i$ induces an isomorphism of graphs between 
$\gamma_i$ and the full subgraph of $\gamma_i^{\infty}$ with vertex set
$\{\phi_i(\m)\}$, where $\m$ runs over all $\Zb/n\Zb$-multisegments.
The description of this last graph follows from \cite{Re97} 
and this finishes the proof of Theorem~\ref{TH1}.
\cqfd

%%%%%%%%%%%%%%%%%%%%%
%
\begin{figure}[t]
\begin{center}
\leavevmode
\epsfxsize =15cm
\epsffile{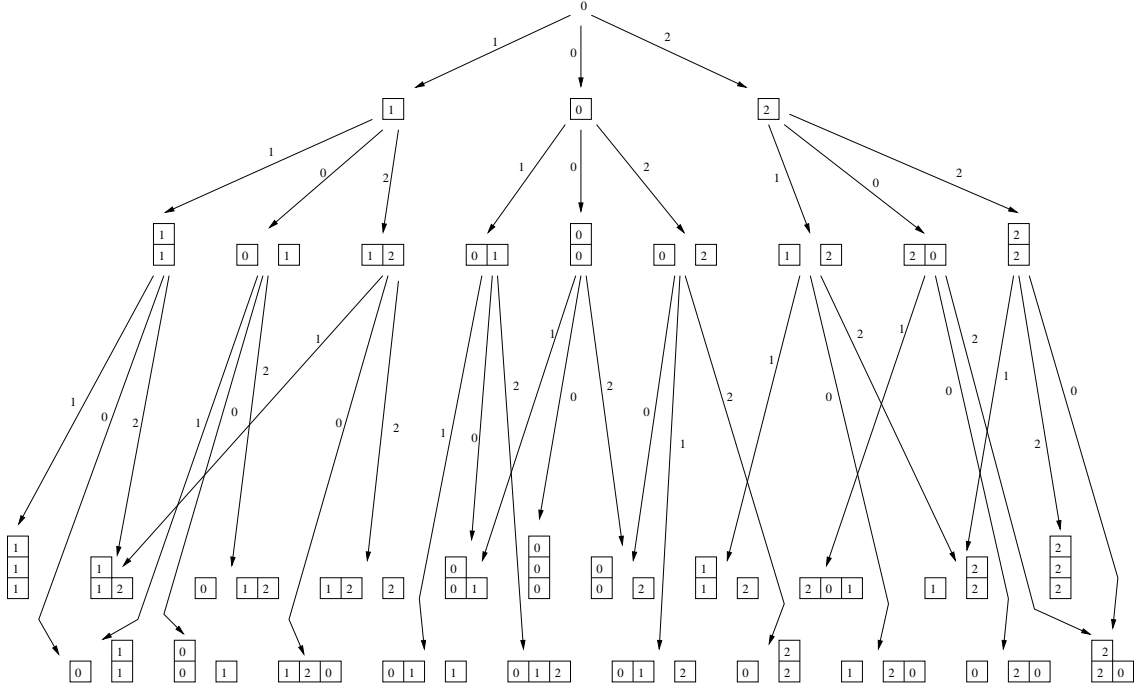}
\end{center}
\caption{\label{FIG2} \small The crystal graph of $U_v^-(\slchap_3)$ 
up to degree 3 in its labelling by aperiodic orbits}
\end{figure}
%
%%%%%%%%%%%%%%%%%%%%%%

\noindent Theorem~\ref{TH1} is illustrated in Figure~\ref{FIG2}
in the case $n=3$.

\bigskip
\noindent
{\bf Remarks 1} \quad It follows easily from Theorem~\ref{TH1} that
the crystal $\widehat{B}(\infty)$ of $\Hc$ decomposes into an infinite number of 
connected components isomorphic to $B(\infty)$, the highest weight vertices
of these components being the periodic multisegments, that is,
the multisegments $\m$ for which $m_{(\ell,r]}=m_{(\ell,s]}$
for all $r,s \in \Zb/n\Zb$.
This is in agreement with \cite{Ka}, Remarks~3.4.10,~3.5.1. 

\smallskip
\noindent
{\bf 2} \quad The crystal graphs of the level $l$ integrable representations
of $\slchap_n$ as described by \cite{JMMO} embed in a natural
way into the crystal graph of $\Hc$ (\cf 4.3 Remark).
Each $l$-tuple of partitions 
$\blambda = (\lambda^1,\ldots , \lambda^l)$ is mapped to
the multisegment obtained by taking the formal sum of the rows
of the $\lambda^i$.

%%%%%%%%%%%%%%%%%%%%%%%%%%%%%%%%%%%%%%%%%%%%%%%%%%%%%%%%%%%%%%%%%%
%  SECTION 5
%%%%%%%%%%%%%%%%%%%%%%%%%%%%%%%%%%%%%%%%%%%%%%%%%%%%%%%%%%%%%%%%%%

\section{Ariki's theorem and Zelevinsky's involution} \label{SECT5}

\subsection{}

Let $K(t)$ denote the graded dual of the $\Nb^{(I)}$-graded space $R(t)$. 
This is made into an associative algebra by taking as product 
the dual of the comultiplication
$$R_k(t)\to \bigoplus_{k=k'+k''} R_{k'}(t)\otimes R_{k''}(t)$$
coming from the restriction maps
$$\H_k(t) \to \H_{k'}(t) \otimes \H_{k''}(t) \, \quad (k=k'+k'').$$
For any $i\in{\bf Z}$ let $\theta_i\in R_1(t)$ be the one-dimensional
module such that $\theta_i(y_1)=t^i$. 
If $t$ is generic (\resp $t=e^{2i\pi/n}$) then
$\{\theta_i\,|\,i\in{\bf Z}\}$ (\resp $\{\theta_i\,|\,i=0,1,...,n-1\}$)
is a basis of $R_1(t)$. The dual basis vectors are denoted by $\theta^i$.
Recall that if $M$ is a simple $\widehat H_k(t)$-module
then $c_k := y_1+y_2+\cdots+y_k$ acts as a scalar, say $z_M$, on $M$.
Then the $i$-restriction of the class $[M]$ is defined in \cite{Ar}
as the class $i$-$\res[M]$ of the $\widehat H_{k-1}(t)$-submodule
$$M_i={\rm Ker}(c_{k-1}-z_M+t^i)^l\subseteq M,\quad l>>1\,,$$
and $i$-$\res$ is extended to $R(t)$ by linearity.
It is straightforward to check that for $M\in R_k(t)$ and 
$f\in R_{k-1}(t)^*$, one has
$
(\theta^i\cdot f)(M)=f(i{-}\res(M)) \,.
$

\begin{lemma} \label{LEM5.1}
{\rm (i)} $\sharp$ induces an automorphism of $K(t)$ such that 
$(\theta^i)^\sharp=\theta^{-i}$.\hfill\break
{\rm(ii)} $\tau$ induces an anti-automorphism of $K(t)$ such that 
$(\theta^i)^\tau=\theta^i$.
\end{lemma}
\Proof
Claim (i) is obvious since $\sharp$ commutes with the natural embedding
$$i:\ \widehat H_{k'}(t)\otimes\widehat H_{k''}(t)\hookrightarrow\widehat H_k(t)\,.$$
Claim (ii) follows from the identity 
$\tau \circ i(a\otimes b)=j(\tau(b)\otimes\tau(a))$
where $a\in\widehat H_{k'}(t)$, $b\in\widehat H_{k''}(t)$, and
$j$ stands for the natural embedding of
$\widehat H_{k''}(t)\otimes\widehat H_{k'}(t)$
into
$\widehat H_k(t)$.
\cqfd

\subsection{}
It was proved by Ariki that the map $f_i\mapsto\theta^i$ 
is an isomorphism of algebras $U^-(\slchap_n)\rightarrow K(t)$ 
if $t=e^{2i\pi/n}$
(\resp $U^-(\Sl_\infty)\rightarrow K(t)$ if $t$ is generic)
(\cite{Ar}, Prop. 4.3; see also the Appendix).
Recall from 2.2 that
$\{[L_O]\}$ is a basis of $R(t)$.
Ariki's theorem states  
that under the previous isomorphism, one has
$b_O \mapsto [L_O]^*$, that is, the canonical basis of $U^-(\slchap_n)$ 
(\resp $U^-(\Sl_\infty)$) 
is mapped to the basis of $K(t)$ dual to $\{[L_O]\}$.

\begin{proposition}\label{PRO5.1}
Let $\m$ be a $\Zb$-multisegment 
{\rm (}\resp an aperiodic \Zb/n\Zb-multisegment{\rm )}
and let $t\in \Cb^*$ be generic {\rm (}\resp $t=e^{2i\pi/n}${\rm )}.
Let $L_{O_\m}=L_\m$ be the corresponding simple $\H_m(t)$-module.
Then, the twisted module $L_\m^\sharp$ is isomorphic to $L_{\m^\sharp}$,
where $\m^\sharp$ is the vertex of the crystal graph of $U_v^-(\Sl_\infty)$
{\rm (}\resp $U_v^-(\slchap_n)${\rm )} obtained from the vertex $\m$ by the 
2-fold symmetry $i \longleftrightarrow -i$ of the graph.
\end{proposition} 
\Proof
Let $\sharp$ denote the automorphism of $U^-(\g)$
induced by the $2$-fold symmetry $i \longleftrightarrow -i$ of
the Dynkin diagram of $\g$.
Repeating the argument of \cite{LLT96} 7.3, 7.4, we see that 
$\sharp$ preserves the canonical basis $\{ b_\m\}$ and maps
$b_\m$ to $b_{\m^\sharp}$, where $\m^\sharp$ is obtained from 
$\m$ by the symmetry $i \longleftrightarrow -i$ of the crystal graph.
The result then follows from Lemma~\ref{LEM5.1} (i) and Ariki's theorem.
\cqfd

%%%%%%%%%%%%%%%%%%%%%
%
\begin{figure}[t]
\begin{center}
\leavevmode
\epsfxsize =10cm
\epsffile{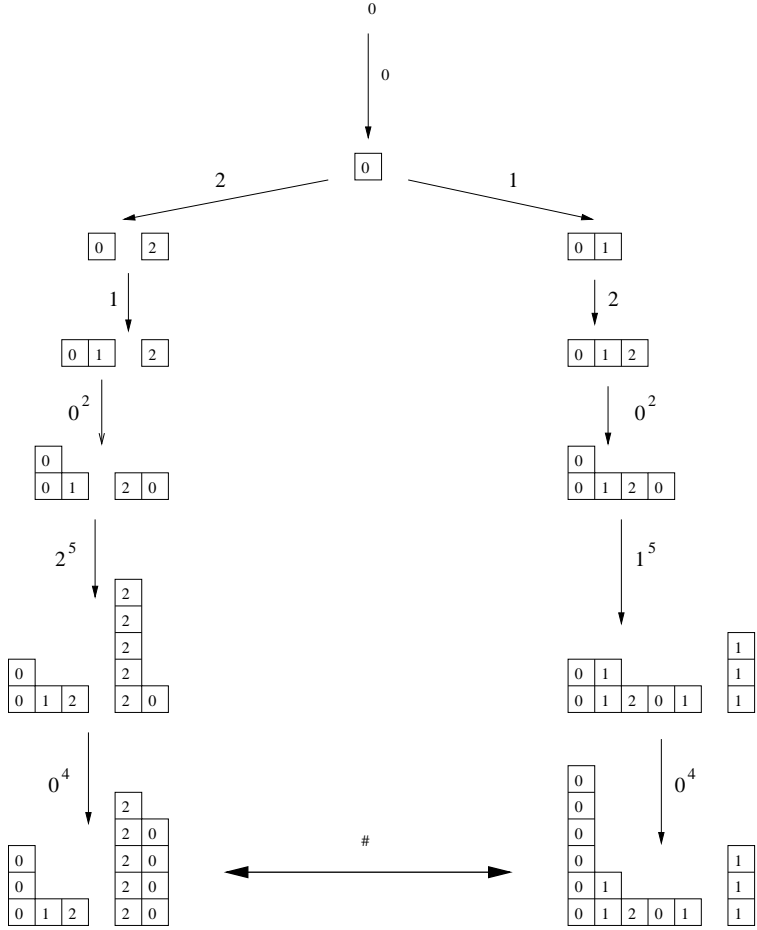}
\end{center}
\caption{\label{FIG3} \small
Computing the dual of a multisegment at a 3rd root of 1}
\end{figure}
%
%%%%%%%%%%%%%%%%%%%%%%

Using Theorem~\ref{TH1}, it is then easy to describe a simple algorithm
for computing $\m^\sharp$. Take any path 
$\emptyset\ {\buildrel i_1\over\longrightarrow} \cdot 
{\buildrel i_2\over\longrightarrow} \cdot \ \cdots \ \cdot
{\buildrel i_k\over\longrightarrow} \ \m$
from the empty multisegment to $\m$ in the crystal graph.
Then $\emptyset\  {\buildrel -i_1\over\longrightarrow} \cdot 
{\buildrel -i_2\over\longrightarrow} \cdot \ \cdots \ 
\cdot {\buildrel -i_k\over\longrightarrow} \ \m^\sharp$
is a path to $\m^\sharp$. This is illustrated in Figure~\ref{FIG3} in the case
$n=3$.

Recall that Zelevinsky's involution $\tau$ is related to $\sharp$ by 
$\tau = \sharp \circ \flat$. As explained in 2.4, $\flat$ acts on a multisegment
$\m$ by simply changing all segments $(\ell;i]$ into $[-i;\ell)$.
Thus Proposition~\ref{PRO5.1} gives as well an algorithm for $\tau$.

\bigskip\noindent
{\bf Remarks 1} \ Consider the crystal of $U^-_v(\Sl_\infty)$.
Among all paths joining a given $\m$ to $\emptyset$, there is a distinguished
one obtained by iterating the following procedure. 
Take the segments $(\ell,i]$ of $\m$ with $i$ minimum, and among them pick
up one with $\ell$ minimum. Then there is always an arrow 
$\m^-_{\,\ell,i}\ {\buildrel i \over\longrightarrow}\ \m$ in the graph.
The algorithm for $\tau$ obtained by using this particular path is
precisely the one described by Moeglin and Waldspurger (\cite{MW}, Lemme II.3). 
Note that here, we need a total order on the vertices $i$ of the quiver,
and there is no obvious analog of this procedure for the cyclic case.

\smallskip\noindent
{\bf 2} \ In \cite{BZ}, Berenstein and Zelevinsky have described for $\g = \Sl_n$
the involutions of the canonical basis of $U^-(\g)$ 
corresponding to $\sharp$ and $\tau$ via Lemma~\ref{LEM5.1} and Ariki's theorem.
They have related them to the transition map between the two indexations
of the basis associated with the two extreme reduced expressions of $w_0$
(\cite{BZ}, Prop. 3.3).
They also noted the coincidence with the multisegment duality coming from
the representation theory of $\Pad$-adic $\GL(m)$ (see \cite{BZ}, 3. Remark).

\subsection{}
Let $\g$ be a symmetrizable Kac-Moody Lie algebra.
We briefly review the geometric construction of $B(\infty)$
conjectured by Lusztig (\cite{Lu90}, 10.4) and 
proved by Kashiwara-Saito \cite{KS}. 
Let $\Gamma$ be the graph associated to the root system
of $\g$. Let $I$ be the set of vertices. Fix an orientation $\Omega$ of $\Gamma$
and set $H=\Omega\cup\overline\Omega$,
where $\overline\Omega$ denotes the orientation opposite to $\Omega$. 
For any arrow $a\in H$ let $\In(a)$ and
$\Out(a)$ be the input and output vertices of $a$. Given a $I$-graded
complex vector space $V=\bigoplus_{i\in I}V_i$, Lusztig has introduced
\cite{Lu90,Lu91}
$$
X_V=\bigoplus_{a\in H}{\rm Hom}(V_{\In(a)},V_{\Out(a)}),
$$
and as above,
$$
E_{V,\Omega}=\bigoplus_{a\in\Omega}{\rm Hom}(V_{\In(a)},V_{\Out(a)}).
$$
The space $E^*_{V,\Omega}$ is naturally identified with $E_{V,\overline\Omega}$,
and, thus, $X_V=E_{V,\Omega} \oplus E_{V,\overline\Omega}$ 
gets identified with the cotangent bundle of $E_{V,\Omega}$.
It is therefore a symplectic manifold, and the natural action of
$G_V$ on $X_V$ is Hamiltonian, with moment map $\mu :\ X_V\rightarrow \g_V$,
where  $\g_V$ is the Lie algebra of $G_V$. 
The $i$-th component of  $\mu$ is the map
$$\mu_i\,:\,X_V\to{\rm End}(V_i),\ \quad B\mapsto
\sum_{a\in H\atop \In(a)=i}\epsilon(a)B_{\overline a}B_a,$$
where $\epsilon(a)=1$ if $a\in\Omega$ and $-1$ if $a\in\overline\Omega$.
The set $\{B \in X_V \,|\, \mu(B) = 0 \}$ is easily seen to be
the union of the conormal bundles $C_O$ of the $G_V$-orbits 
$O \in E_{V,\Omega}$.  

An element $B\in X_V$ is nilpotent if there is an integer $l$ such that
$B_{a_1}\circ B_{a_2}\circ\cdots\circ B_{a_l}=0$ for all
sequences $(a_1, a_2,...,a_l)$ such that the above expression is meaningful.
Define
$$
\Lambda_V=\{B\in X_V\,|\,\mu(B)=0\mbox{\ and\ }B\ \mbox{is\ nilpotent}\}.
$$
It is known \cite{Lu90,Lu91} that $\Lambda_V$ is a closed equidimensional 
$G_V$-stable subvariety of $X_V$. 
For $\Gamma$ of type $A_{n-1}$, $A_\infty$ or $A_{n-1}^{(1)}$, an element $B$ of
$\Lambda_V$ can be identified with a pair $(x,y)$ of commuting nilpotent 
endomorphisms of $V$ of degree $-1$ and $+1$ respectively.
Moreover, the irreducible components of $\Lambda_V$ are the 
closures $\overline{C_O}$ of the conormal bundles to the nilpotent $G_V$-orbits 
$O$ in $E_{V,\Omega}$
for type $A_{n-1}$ or $A_\infty$ (\cite{Lu91}, 14), and to the aperiodic
nilpotent $G_V$-orbits for type $A_{n-1}^{(1)}$ (\cite{Lu91}, 15).

To the $I$-graded space $V$ we associate the negative root 
$\nu=-\sum_i({\rm dim}V_i)\alpha_i.$ 
%Put $X_\nu=X_V$,
%$E_{\nu,\Omega}=E_{V,\Omega}$, $\Lambda_\nu=\Lambda_V$, $G_\nu=G_V$.
Let $B(\infty;\nu)$ be the set of irreducible components of $\Lambda_V$. 
Lusztig has introduced a crystal structure on the set $\sqcup_\nu B(\infty;\nu)$,
and conjectured that this graph is isomorphic to the crystal $B(\infty)$
of $U^-_v(\g)$ (\cite{Lu90}, 8.9). 
This was proved by Kashiwara-Saito \cite{KS}. More precisely, 
recall that the vertices of $B(\infty)$ are identified with orbits $O \in {\cal O}$
for type $A_{n-1}$ or $A_\infty$, and aperiodic orbits $O \in {\cal O}$ 
for type $A_{n-1}^{(1)}$.
Then the map $\overline{C_O} \mapsto O$ is an isomorphism of crystals
$\sqcup_\nu B(\infty;\nu) \rightarrow B(\infty)$.

The algebra $U_v(\g)$ admits an anti-automorphism $*$
such that $e_i^*=e_i$, $f_i^*=f_i$, and $(v^h)^*=v^{-h}$,
which induces an involution (still denoted by $*$) of $B(\infty)$ \cite{Ka}.
It is shown in \cite{KS}, Section 5.3 that
the involution of $B(\infty;\nu)$ induced by the transpose map $B\mapsto {}^tB$ 
on $\Lambda_V$ coincides with $*$.
(This was also conjectured by Lusztig \cite{Lu90}, 10.4.)
More precisely,
if $x\in O\subseteq N_{V,\Omega}$ then $*(\overline{C_O})$ is the irreducible
component containing $({}^ty,{}^tx)$, where $y\in N_{V,\overline\Omega}$ is
generic such that $[x,y]=0$. Summarizing the discussion,
and taking into account Lemma~\ref{LEM5.1} (ii) and Ariki's theorem, 
we can state:

\begin{proposition}
Let $\Gamma$ be of type $A_{n-1}^{(1)}$ or $A_\infty$.
For $x$ in $E_{V,\Omega}$,
set 
$$
Z(x)=\{y\in E_{V,\overline\Omega}\ |\ x\circ y=y \circ x\}\,.
$$
Assume that $x$ is nilpotent (and aperiodic in type $A_{n-1}^{(1)}$).
Then, \hfill\break
{\rm (i)} there is a unique $G_V$-orbit $O'\subset N_{V,\overline\Omega}$ 
such that $O'\cap Z(x)$ is Zarisky dense and open in $Z(x)$, and this orbit 
only depends on the orbit $O$ of $x$;\hfill\break
{\rm (ii)} if $O=O_{\bf m}$ and $y \in O'$, then  ${}^ty \in O_{\m^\tau}$,
where the multisegment ${\bf m}^\tau$ is obtained from $\m$ as explained in 
{\rm 5.2}.
\end{proposition}

For $\Gamma$ of type $A_\infty$, we recover the description of
Zelevinsky's involution conjectured in \cite{Zel81} and established
in \cite{MW}. 

%%%%%%%%%%%%%%%%%%%%%%%%%%%%%%%%%%%%%%%%%%%%%%%%%%%%%%%%%%%%%%
% APPENDIX
%%%%%%%%%%%%%%%%%%%%%%%%%%%%%%%%%%%%%%%%%%%%%%%%%%%%%%%%%%%%%%

%%%%%%%%%%%%%%%%%%%%%%%%%%%%%%%%%%%%%%%%%%%%%%%%%%%%%%%%%%%%%%%%%%%%%%%%%%%%%

\section{Appendix}
In \cite{Ar} Ariki uses the non-vanishing theorem
for affine Hecke algebras at roots of unity announced in \cite{Groj}. The purpose
of this appendix is to indicate a simple proof of this non-vanishing 
result for type $A$. S. Ariki has informed one of us (E.V.) 
that this simple argument was known to him and to G. Lusztig. 
For the sake of completeness, we include it here.

Fix $m\in\Nb^*$ and put $t=e^{2i\pi/n}$.
We use the same notations as in Sections 3.1 and 5.3. In particular
$(\Gamma,\Omega)$ is the cyclic quiver with $n$ vertices and 
$V$ is a $I$-graded vector space of dimension ${\bf d}\in\Nb^{(I)}$.
Let $S_{\bf d}$ be the set of sequences ${\bf i}=(i_1,i_2,...,i_m)\in I^m$ 
such that $\sharp\{k\,|\,i_k=i\}=d_i$. Given such a sequence
let $F_{\bf i}$ be the variety of flags of type ${\bf i}$ in $V$ 
(see \cite{Lu91}, Section 1). Let $\tilde F_{\bf i}$ be the variety of pairs 
$(\phi,x)\in F_{\bf i}\times E_{V,\Omega}$ such that $\phi$ is $x$-stable.
As usual, we denote by ${\cal D}(X)$ the bounded derived category of 
complexes of $\Cb$-sheaves on a complex algebraic variety $X$.
Let $\pi_{\bf i}\,:\,\tilde F_{\bf i}\to E_{V,\Omega}$ be the second
projection and put 
${\cal L}_{\bf i}=(\pi_{\bf i})_!(\Cb)\in{\cal D}(E_{V,\Omega})$.
By definition the degree-${\bf d}$ component of
Lusztig's canonical basis $\{b_O\}$ is labelled by the
$G_V$-orbits $O\subset N_{V,\Omega}$ such that a shift of
the intersection cohomology complex $IC_O$ is a direct factor in 
$\bigoplus_{{\bf i}\in S_{\bf d}}{\cal L}_{\bf i}$ (use \cite{Lu91}, Section 2.4).

Let $F$ be the variety of complete flags in $\Cb^m$. 
The cotangent variety $T^*F$ is identified with the set of pairs 
$(\phi,x)\in F\times{\rm End}(\Cb^m)$ such that $\phi$ is $x$-stable. 
If $s\in GL(m)$ is semi-simple let $G_s\subseteq GL(m)$
be its centralizer, let $N_s\subset{\rm End}(\Cb^m)$ 
be the set of nilpotent matrices $x$ such that
$sxs^{-1}=t^{-1}x$ and let $\tilde F_s$ be the set of pairs $(\phi,x)\in T^*F$
such that $x\in N_s$ and the flag $\phi$ is fixed by $s$. 
Let $\pi_s\,:\,\tilde F_s\to N_s$ be the second projection
and set ${\cal L}_s=(\pi_s)_!(\Cb)\in{\cal D}(N_s)$.
Given a sequence ${\bf d}\in\Nb^{(I)}$ such that $\sum_id_i=m$, if
$s$ is a diagonal matrix such that $t^i$ has the multiplicity $d_i$ in
the spectrum of $s$ then $G_s$, $N_s$ are identified with $G_V$, $N_{V,\Omega}$
above. In this case set ${\cal L}_{\bf d}={\cal L}_s$.
Then, Ginzburg's geometric construction of the simple modules of the
affine Hecke algebra (see \cite{CG}) implies that the simple modules of 
$\widehat H_m(t)$ are labelled by the $G_V$-orbits $O\subset N_{V,\Omega}$
such that a shift of the complex $IC_O$ is a direct factor in
$\bigoplus_{\bf d}{\cal L}_{\bf d}$.

In our case the non-vanishing result is precisely that the orbits labelling
the degree-$m$ component of the basis $\{b_O\}$ are the same as
the orbits labelling the simple $\widehat H_m(t)$-modules. It is a 
consequence of the obvious equality of the complexes ${\cal L}_{\bf d}$ and
$\bigoplus_{{\bf i}\in S_{\bf d}}{\cal L}_{\bf i}$
(simply observe that $F_{\bf d} = \sqcup_{{\bf i} \in S_{\bf d}} F_{\bf i}$).

%%%%%%%%%%%%%%%%%%%%%%%%%%%%%%%%%%%%%%%%%%%%%%%%%%%%%%%%%%%%%%
\bigskip\noindent
{\bf Acknowledgements} \quad We are grateful to M.-F. Vigneras for posing
the question which motivated this work.
B.L. and J.-Y.T. would like to thank A. Zelevinsky for his beautiful
series of lectures on ``Representations of quantum groups and
piecewise linear combinatorics'' given at the Universit\'e de Marne-la-Vall\'ee
in 1994, which arose their interest in this field.

E.V. is partially supported by EEC grant no ERB FMRX-CT97-0100.

%%%%%%%%%%%%%%%%%%%%%%%%%%%%%%%%%%%%%%%%%%%%%%%%%%%%%%%%%%%%%%%%%
% BIBLIOGRAPHY
%%%%%%%%%%%%%%%%%%%%%%%%%%%%%%%%%%%%%%%%%%%%%%%%%%%%%%%%%%%%%%%%%

%\bigskip\bigskip

%\bigskip\bigskip

\end{document}